\documentclass[12pt,reqno]{amsart}

\usepackage{color}
\usepackage{amsmath, amsthm, amssymb}
\usepackage{amsfonts}
\usepackage[ansinew]{inputenc}
\usepackage[dvips]{epsfig}
\usepackage{graphicx}
\usepackage[english]{babel}
\usepackage{hyperref}

\usepackage{cite}
\usepackage{graphicx}
\usepackage{amscd}
\usepackage{color}
\usepackage{bm}
\usepackage{enumerate}

\usepackage{verbatim}
\usepackage{hyperref}
\usepackage{amstext}
\usepackage{latexsym}
\theoremstyle{plain}
\newtheorem{Theorem}{Theorem}[section]

\newtheorem{Corollary}[Theorem]{Corollary}
\newtheorem{Proposition}[Theorem]{Proposition}
\newtheorem{Lemma}[Theorem]{Lemma}
\newtheorem{Remark}[Theorem]{Remark}

\numberwithin{Theorem}{section}
\numberwithin{equation}{section}

\def\proof{\noindent{{\bf Proof. }}}
\def\square{\vbox{
\hrule height .4pt \hbox{\vrule width .4pt height 7pt \kern 7pt
\vrule width .4pt} \hrule height .4pt }}
\def\QED{\hfill {$\square$}\goodbreak \medskip}

\newcommand{\average}{{\mathchoice {\kern1ex\vcenter{\hrule height.4pt
width 6pt depth0pt} \kern-9.7pt} {\kern1ex\vcenter{\hrule
height.4pt width 4.3pt depth0pt} \kern-7pt} {} {} }}

\def\R{\mathbb{R}}
\def\loc{\text{loc}}

\renewcommand{\a }{\alpha }

\newcommand{\e }{\varepsilon }

\renewcommand{\l }{\lambda }

\newcommand{\vp }{\varphi }

\newcommand{\s }{\sigma }

\renewcommand{\t }{\tau }

\renewcommand{\O }{\Omega }

\newcommand{\ov}{\overline}

\newcommand{\be}{\begin{equation}}
\newcommand{\ee}{\end{equation}}

\newcommand{\de}{\partial}

\DeclareMathOperator{\inv}{inv}

\newcommand{\N}{\mathbb{N}}

\newcommand{\Z}{\mathbb{Z}}

\newcommand{\cG}{{\mathcal G}}
\newcommand{\cH}{{\mathcal H}}
\newcommand{\cI}{{\mathcal I}}

\newcommand{\cL}{{\mathcal L}}

\newcommand{\cO}{{\mathcal O}}
\newcommand{\cP}{{\mathcal P}}

\newcommand{\cU}{{\mathcal U}}

\newcommand{\eps}{\varepsilon}

\DeclareMathOperator{\id}{id}

\renewcommand{\epsilon}{\varepsilon}

\begin{document}

\title[Periodic Serrin]
{Unbounded periodic solutions to Serrin's overdetermined boundary value problem}

\author{Mouhamed Moustapha Fall}
\address{M. M. F.: African Institute for Mathematical Sciences in Senegal, KM 2, Route de
Joal, B.P. 14 18. Mbour, Senegal.}
\email{mouhamed.m.fall@aims-senegal.org}

\author{Ignace Aristide Minlend}
\address{I. A. M.:African Institute for Mathematical Sciences in Senegal, KM 2, Route de
Joal, B.P. 14 18. Mbour, Senegal. }
\email{ignace.a.minlend@aims-senegal.org}

\author{Tobias Weth}
\address{T.W.:  Goethe-Universit\"{a}t Frankfurt, Institut f\"{u}r Mathematik.
Robert-Mayer-Str. 10 D-60054 Frankfurt, Germany.}

\email{weth@math.uni-frankfurt.de}



\begin{abstract}
We study the existence of nontrivial unbounded domains $\Omega$ in $\R^N$ such that the overdetermined problem
$$
-\Delta u = 1 \quad \text{in $\Omega$}, \qquad u=0, \quad \partial_\nu u=\textrm{const}  \qquad \text{on $\partial \Omega$}
$$
admits a solution $u$. By this, we complement Serrin's classification result from 1971 which
yields that every bounded domain admitting a solution of the above problem is a ball in $\R^N$. The domains we construct are periodic in some variables and radial in the other variables, and they bifurcate from a straight (generalized) cylinder or slab. We also show that these domains are uniquely self Cheeger relative to a period cell for the problem.
\end{abstract}

\maketitle
\section{Introduction and main result}
In 1971, Serrin \cite{Serrin}
established a celebrated result on the overdetermined problem of finding a domain $\O \subset \R^N$ and a $C^2$-function $u: \overline \Omega \to \R$ such that
\begin{equation}\label{eq:Pro-1}
 -\Delta u=1 \qquad \text{ in $\Omega$}
\end{equation}
and
\begin{equation}\label{eq:Pro-2}
u=0, \quad \partial_\nu u=\textrm{const}  \qquad \text{on $\partial \Omega$.}
\end{equation}
Here  $\nu$ is the unit outer normal on  $\partial \O$. More precisely, in \cite{Serrin} Serrin proved that, if $\O$ is a bounded domain of class $C^2$ such that (\ref{eq:Pro-1}), (\ref{eq:Pro-2}) admits a solution,  then $\O$ is a ball. The problem (\ref{eq:Pro-1}),~(\ref{eq:Pro-2}) arises in e.g.  in fluid dynamics and the linear theory of torsion, and we refer the reader to  \cite{Sirakov,Serrin} for a detailed account on its relevance. As we shall discuss further below, it is also related to the notion of Cheeger sets, which in turn has  applications in the denoising  problem in  image processing.  The proof of Serrin's classification result for (\ref{eq:Pro-1}),~(\ref{eq:Pro-2}) relies on the moving plane method, and it extends to the much more general problem where (\ref{eq:Pro-1}) is replaced by
\begin{equation}\label{eq:Pro-1-general}
 -\Delta u=f(u), \quad u>0 \qquad \text{ in }\quad \Omega
\end{equation}
with a locally Lipschitz continuous nonlinearity $f$. We note that the positivity assumption in (\ref{eq:Pro-1-general}) is essential, and  by the strong maximum principle it is automatically satisfied for solutions of (\ref{eq:Pro-1}),~(\ref{eq:Pro-2}). The moving plane method,  which Serrin established in a PDE context, is inspired by Alexandrov's reflection principle \cite{Alexandrov} for constant mean curvature hypersurfaces. On the other hand, Weinberger \cite{Wein} found  a simpler argument to prove Serrin's result for problem (\ref{eq:Pro-1}),~(\ref{eq:Pro-2}) without the moving plane method, but his argument does not cover the more general problem ~(\ref{eq:Pro-2}),~(\ref{eq:Pro-1-general}).

The result of Serrin parallels an earlier important result by Alexandrov \cite{Alexandrov} which states that closed embedded hypersurfaces with constant mean curvature (CMC hypersurfaces in short) are round spheres. This rigidity result for bounded embedded CMC hypersurfaces stands in striking contrast to the rich structure of unbounded CMC hypersurfaces which has been explored
in classical and more recent papers. For a survey, we refer the reader to \cite{mazzeo}.
We recall in particular that already in 1841, Delaunay \cite{Delau} constructed and classified unbounded surfaces of revolution in $\R^3$ with constant mean curvature. As Delaunay's construction shows, these surfaces bifurcate from a straight cylinder (see also \cite[Section 2]{ScSi} for a different proof of the latter statement).

For some time, it was unknown whether the problem (\ref{eq:Pro-2}),~(\ref{eq:Pro-1-general}) admits solutions in nontrivial unbounded domains. In fact, in \cite{BCN}, Berestycki, Caffarelli and Nirenberg conjectured that if $\Omega \subset \R^N$ is an {\em unbounded} sufficiently regular
domain such that $\R^N \setminus \overline \Omega$ is connected and $f: [0,\infty) \to \R$ is a local Lipschitz continuous function such that the overdetermined problem (\ref{eq:Pro-2}),~(\ref{eq:Pro-1-general}) admits a solution, then $\Omega$ is an affine half space or the complement $B^c$ of a ball $B \subset \R^N$ or a product of the form $\R^j \times B$ bzw. $\R^j \times B^c$  with a ball $B \subset \R^{N-j}$ (after a suitable rotation).

This conjecture has been disproved by Sicbaldi \cite{Sic} in dimensions $N \ge 3$. More precisely, in the case where $f(u)= \lambda_1 u$ with $\lambda_1>0$ suitably chosen, it was shown in \cite{Sic}
that there exist periodic domains of revolution such that the problem (\ref{eq:Pro-1-general}),~(\ref{eq:Pro-2}) admits a positive solution. Moreover, these domains bifurcate from the straight cylinder $\R \times B$, where $B \subset \R^{N-1}$ is a ball. The construction in \cite{Sic} relies on topological degree theory and therefore does not give rise to a smooth branch of domains; moreover, the case $N=2$ was not included. Later in \cite{ScSi}, Sicbaldi and Schlenk extended the result to dimensions $N \ge 2$, and they obtained a smooth branch of domains via the use of the Crandall-Rabinowitz bifurcation theorem (see \cite{M.CR}).

In further recent papers, different types of nonlinearities $f$ and unbounded domains have been
considered in the context of the general Serrin problem  (\ref{eq:Pro-2}),~(\ref{eq:Pro-1-general}). In \cite{hauswirth-et-al,Traizet}, the authors study examples of unbounded domains where (\ref{eq:Pro-2}),~(\ref{eq:Pro-1-general})  is solvable with $f= 0$, i.e. with harmonic functions. Moreover, in \cite{del-Pino-pacard-wei},  the authors consider a monostable nonlinearity $f$, and they construct domains whose boundary is close to dilations of a given CMC-hypersurface and such that  (\ref{eq:Pro-2}),~(\ref{eq:Pro-1-general}) is solvable. It is clear from these works that the existence and shape of such domains depend in a crucial way on the function $f$. For negative results, excluding the solvability of (\ref{eq:Pro-2}),~(\ref{eq:Pro-1-general}) in nontrivial unbounded domains belonging to certain domain classes (e.g. epigraphs), we refer the reader to \cite{farina-valdinoci,farina-valdinoci:2010-1,farina-valdinoci:2010-2,farina-valdinoci:2013-1,Ros-Sicbaldi,Ros-Ruiz-Sicbaldi} and the references therein.

In the present paper, we wish to analyze the original form of Serrin's problem (\ref{eq:Pro-1}),~(\ref{eq:Pro-2}), i.e. the case $f \equiv 1$, in unbounded domains. More precisely, we study domains of the form
$$
\O_\phi:=\left\{ (z,t)\in
\mathbb{R}^{n}\times\mathbb{R}^m\,:\, |z|<\phi(t)  \right\}\subset\R^N,
$$
where $N=n+m$ and $\phi: \R^m \to (0,\infty)$ is an even and $2\pi \Z^m$-periodic function. The following is our main result.

\begin{Theorem}\label{teo1}
For each $n,m \geq 1$ and $\alpha \in (0,1)$, there exists $\lambda_{*}=\lambda_*(n)>0$ and a smooth map
\begin{align*}
 (-\e_0,\e_0) &\longrightarrow  (0,\infty) \times C^{2,\alpha}(\mathbb{R}^m)\\
        s&\longmapsto (\lambda_s,\vp_s)
 \end{align*}
with $\vp_0\equiv 0$, $\lambda_0=\lambda_*$ and such that
for all
$s\in(-\e_0,\e_0),$ letting $\phi_s=\l_s+\vp_s$,  there exists a solution $u \in C^{2,\a}(\ov{\O_{\phi_s}})$ of the overdetermined problem
\begin{equation}
 \label{eq:Proe1}
 \left\{
   \begin{aligned}
  -\Delta u&=1&&\qquad \textrm{in}\quad \Omega_{\phi_s}\\
u&=0,\quad  \partial_\nu u=-\frac{\lambda_s}{n} &&\qquad \textrm{on }\quad \partial\Omega_{\phi_s}
 \end{aligned}
\right.
\end{equation}
 in the domain
\begin{equation}\label{petcy}
\O_{\phi_s}= \biggl\{(z,t)\in \mathbb{R}^{n}\times\mathbb{R}^m:\quad
|z|<  \phi_s( t)\biggl\}.
\end{equation}
Moreover,  for every $s \in (-\e_0,\e_0)$, the function $\vp_s$ is even in $t_1,\dots, t_m$, ${2\pi} $-periodic in
$t_1,\dots, t_m$ and invariant with respect to permutations of the
variables $t_1,\dots,t_m$. Furthermore, we
have
$$
\vp_s(t)=s\Bigl(\sum_{j=1}^m \cos(t_j)+\mu_s(t)\Bigr) \qquad \text{for $s \in (-\e_0,\e_0)$}
$$
with a smooth map $ (-\e_0,\e_0) \to  C^{2,\alpha}(\mathbb{R}^m)$, $s\mapsto \mu_s$ satisfying
$$
\int_{[0,2\pi]^m} \mu_s(t) \cos( t_j) \,dt = 0 \qquad \text{for $s \in (-\e_0,\e_0)$, $j=1,\dots,m$,}
$$
and $\mu_0\equiv 0$.
\end{Theorem}

Note here that, since the domain $\O_{\phi_s}$ is radially symmetric in $z$ for every fixed $s \in (-\e_0,\e_0)$,
the corresponding solution $u$ is also radially symmetric in the $z$-variable.
We also remark that the bifurcation value $\lambda_*$ in Theorem~\ref{teo1} is the unique  zero of the function
$$
\sigma: (0,\infty) \to \R,\quad
\sigma(\rho) =
\left\{
  \begin{aligned}
 &\rho \tanh (\rho)-1,&&\qquad \text{in case $n=1$,}\\
 &\frac{1}{n}\Bigl(\rho \frac{I_{\nu+1}(\rho)}{I_\nu(\rho) }-1 \Bigr),&&\qquad
\text{in case $n \ge 2$,}\\
  \end{aligned}
\right.
$$
see Propositions~\ref{proposition-eigenvalues} and~\ref{propPhi} below.
 Here $I_\nu$ is the modified Bessel function of the first kind of order $\nu= \frac{n-2}{2}$.
Numerically, $\lambda_*$ is given by
$$
\lambda_*  \approx
\left\{
  \begin{aligned}
 &1.199 &&\qquad \text{in case $n=1$;}\\
 &1.608&&\qquad
\text{in case $n = 2$;}\\
 &1.915&&\qquad
\text{in case $n = 3$.}
  \end{aligned}
\right.
$$
As remarked above, Sicbaldi and Schlenk \cite{ScSi} have derived --
in the special case $m=1$ -- a result analogous to
Theorem~\ref{teo1}  where (\ref{eq:Pro-1}) is replaced by
(\ref{eq:Pro-1-general}) with $f(u)=\lambda_1 u$. Our proof
of Theorem~\ref{teo1} is partly inspired by \cite{ScSi} and also relies on
the Crandall-Rabinowitz Theorem, but there are key differences due to
the special form of (\ref{eq:Pro-1}). We believe that our approach
can also be generalized to study Serrin's overdetermined problem on
Riemannian manifolds. Related to this, we mention the recent works \cite{MS,farina-valdinoci:2013-2}. In \cite[Theorem 5]{farina-valdinoci:2013-2}, necessary conditions for the solvability of some overdetermined problems on manifolds are given, and  in \cite{MS} the case $f(u)= \lambda_1 u$ is studied in the product
manifolds $S^N \times \R$ and $H^N \times \R$. 

 The overdetermined problem (\ref{eq:Pro-1}),~(\ref{eq:Pro-2}) is related to a generalized notion of Cheeger sets. To define this generalized notion, let $S$, $\Omega$ be open subsets of $\R^N$. For a subset $A \subset S$ with Lipschitz boundary, we let $P(A,S):= H^{N-1}(\de A \cap S)$ denote the relative perimeter of $A$ in $S$. Here and in the following, $H^{N-1}$ denotes the $N-1$-dimensional Hausdorff measure. For an equivalent definition which extends to Borel subsets $A$ of $S$, see e.g. \cite[Def. 13.6]{leoni}.  We then define the
Cheeger constant of $\Omega$ relative to $S$ as \be
\label{eq:def_Cheeg} h(\O,S):= \inf_{A\subset \O \cap
S}\frac{P(A,S)}{|A|}, \ee where the infimum is taken over subsets $A
\subset \Omega \cap S$ with Lipschitz boundary. If this constant is
attained by some subset $A \subset \O \cap S$ with Lipschitz
boundary, then $A$ will be called a Cheeger set of $\O$ relative to
$S$. If $\Omega$ has a Lipschitz boundary and $A=\Omega  \cap S$
attains the constant $h(\O,S)$ in \eqref{eq:def_Cheeg},  we say that
$\O$ is self-Cheeger relative to $S$. Moreover, if $A=\O \cap S$ is
the only set which attains $h(\Omega,S)$, we say that $\O$ is
uniquely self-Cheeger relative to $S$. These notions generalize the
classical notions of the Cheeger constant and Cheeger sets which
correspond to the case $S= \R^N$, see e.g. \cite{Leon}.

We have the following corollary of Theorem~\ref{teo1}.

\begin{Corollary}\label{cor:Cheeger}
For every $s\in (-\e_0,\e_0)$ and $a, b \in \pi \Z^m$ with $a_i < b_i$ for $i=1,\dots,m$, the set $\O_{\phi_s}$ given in Theorem~\ref{teo1} is uniquely self-Cheeger  relative to the set
\begin{equation}
  \label{eq:defsalphabeta}
 S_a^b:= \R^n \times (a_1, b_1) \times \dots \times (a_m, b_m)  \quad \subset \quad \R^N= \R^n \times \R^m
\end{equation}
with corresponding relative Cheeger constant $h(\Omega_{\phi_s},S_{a}^b)=
\dfrac{n}{\lambda_s}.$
\end{Corollary}
 As discussed
in detail in the illuminating surveys \cite{Leon,EParini}, self-Cheeger sets arise in various problems as e.g. the
construction of prescribed mean curvature graphs or the
regularization of noisy images within the ROF model.

The link between Serrin's over-determined problem  (\ref{eq:Pro-1}),~(\ref{eq:Pro-2}) and Cheeger sets
on $N$-dimensional Riemannian manifolds was also studied by the second author in
\cite{I.A.M}, where he proved the existence of a family of uniquely
self-Cheeger sets $(\O_\e)_{\e\in (0,\e_0)}$ with classical Cheeger
constant $h(\O_\e)=\frac{N}{\e}$.

The paper is organized as follows. In Section~\ref{sec:transf-probl-its} we transform the overdetermined problem to an equivalent boundary value  problem on a fixed underlying domain with a $\phi$-dependent metric. In Section~\ref{proposition-eigenvalues}, we then study the eigenvalues and eigenfunctions of the linearization of the problem at constant functions $\phi \equiv \lambda$. In particular, we study the dependence of the eigenvalues on $\lambda>0$. In Section~\ref{sec:compl-proof-theor}, we then complete the proof of Theorem~\ref{teo1} via the Crandall-Rabinowitz Theorem. Finally, in Section~\ref{sec:Per}, we give the proof of Corollary~\ref{cor:Cheeger}.

\section{The transformed problem and its linearization}
\label{sec:transf-probl-its}
We fix $\alpha \in (0,1)$ in the following.  For $j \in \N \cup \{0\}$,  we consider the Banach space
$$
C_{p,e}^{j,\alpha}(\mathbb{R}^m):=\biggl\{  \phi \in C^{j,\alpha}(\mathbb{R}^m) \::\: \textrm{$\phi$ is even and $2\pi$-periodic in $t_1,\dots,t_m$}\biggl\}
$$
Let
$$
\cU:= \{ \phi \in C^{2,\alpha}_{p,e}(\R^m) \::\: \phi>0 \}.
$$
For a function $\phi \in \cU$, we define
$$
\O_\phi:=\left\{ (z,t)\in
\mathbb{R}^{n}\times\mathbb{R}^m\,:\, |z|<\phi(t)  \right\}
$$
as well as the spaces
\begin{align*}
&C_{p,e}^{j,\alpha}(\O_\phi, \R^{k}):=\bigl\{  u \in C^{j,\alpha}(\Omega_\phi,\R^{k}) : \textrm{$u$ is even and $2\pi$-periodic in $t_1,\dots,t_m$}\bigl\},\\
&C_{p,e}^{j,\alpha}(\partial \O_\phi, \R^{k}):=\bigl\{  u \in
C^{j,\alpha}(\partial \Omega_\phi,\R^k) : \textrm{$u$ is even and
$2\pi$-periodic in $t_1,\dots,t_m$}\bigl\}
\end{align*}
for $j=0,1,2$, $k \in \N$. If $k=1$, we simply write $C_{p,e}^{j,\alpha}(\overline \O_\phi)$ and $C_{p,e}^{j,\alpha}(\partial \O_\phi)$. Moreover, in the special case $\phi \equiv 1$ we write
$$
\O:= \O_1= \left\{ (z,t)\in
\mathbb{R}^{n}\times\mathbb{R}^m\,:\, |z|<1  \right\}.
$$
Every $\phi \in \cU$ gives rise to a locally $C^{2,\alpha}$-regular map
\be \label{eq:def-diffeom-Psi_phi}
\Psi_\phi: \R^{n+m} \to \R^{n+m}, \qquad \Psi_\phi(z,t)= (\phi(t)z,t).
\ee
such that $\Psi_\phi$ maps $\O$ diffeomorphically onto $\O_\phi$. Let the metric $g_\phi$ be
 defined as the pull back of the euclidean metric $g_{eucl}$ under the map $\Psi_\phi$, so
 that $\Psi_{\phi}:(\overline \O,g_\phi) \to (\overline \O_{\phi},g_{eucl})$ is an isometry.
  Hence  our original problem is equivalent to the overdetermined problem consisting of the Dirichlet problem
\begin{align}\label{eq:Proe1-1}
  \begin{cases}
    -\Delta_{g_\phi} u= 1& \quad \textrm{ in}\quad\O\\
u=0&  \quad\textrm{ on} \quad \partial \O\\
  \end{cases}
  \end{align}
and the additional Neumann condition
\begin{equation}
\label{eq:Proe1-1-neumann}
\partial_{\nu_{\phi}} u \equiv -c \quad\textrm{ on }\quad \partial \O.
\end{equation}
Here
$$
\nu_\phi: \partial \O \to  \R^{n+m}
$$
is the unit outer normal vector field on $\partial \O$ with respect to $g_\phi$. Since \mbox{$\Psi_{\phi}: (\overline \O, g_\phi) \to (\overline \O_\phi, g_{eucl})$} is an isometry, we have
\begin{equation}
  \label{eq:rel-mu-phi-nu-phi}
\nu_{\phi} = [d \Psi_\phi]^{-1} \mu_\phi \circ \Psi_\phi  \qquad \text{on $\partial \Omega$,}
\end{equation}
where $\mu_\phi: \partial \O_\phi \to \R^{n+m}$ denotes the outer normal on $\partial \Omega_\phi$ with respect to the euclidean metric $g_{eucl}$ given by
\begin{equation}
  \label{eq:def-mu-phi}
\mu_\phi(z,t) = \frac{(\frac{z}{|z|}, -\nabla \phi(t))}{\sqrt{1+|\nabla \phi(t)|^2}} \in \R^{n+m} \qquad \text{for $(z,t) \in \partial \Omega_\phi.$}
\end{equation}
Here and in the following, we distinguish different types of derivatives in our notation. If $f: O \to \R^\ell$ is a $C^1$-map defined on an open set $O \subset \R^k$, we write $d f(x) \in \cL(\R^k, \R^\ell)$ for its derivative at a point  $x \in O$. In contrast, we shall use the symbols $D$ or $D_\phi$ to denote functional derivatives. More precisely, if $X,Y$ are infinite dimensional normed (function) spaces and $F \in C^1(O, Y)$, where $O \subset X$ is open, we let $D F(\phi)$ or $D_\phi F(\phi) \in \cL(X,Y)$ denote the Fr\'echet derivative of $F$ at a function $\phi \in O$.

 The following lemma is concerned with the well-posedness of problem (\ref{eq:Proe1-1}).
\begin{Lemma}
\label{sec:peri-solut-serr-1}
For any $\phi \in \cU$, there is a unique solution $u_\phi \in C^{2,\alpha}_{p,e}(\overline \O)$ of (\ref{eq:Proe1-1}), and the map
\begin{equation}
  \label{eq:def-smooth-map}
C^{2,\alpha}_{p,e}(\R^m) \to C^{2,\alpha}_{p,e}(\overline \O),\qquad
\phi \mapsto u_\phi
\end{equation}
is smooth. Moreover we have the following properties.
\begin{itemize}
\item[(i)] For any $\phi \in \cU$, the functions $u_\phi: \overline \Omega \to \R$ and $\partial_{\nu_\phi} u_\phi: \partial \O \to \R$
are radially symmetric in the $z$-variable.
\item[(ii)] For a constant function $\phi \equiv \lambda>0$, we have $u_\lambda(z,t)= \frac{\lambda^2-|\lambda z|^2}{2n}$.
\item[(iii)] Let $\cP \subset \cL(\R^m)$ denote the subset of all coordinate permutations in $t_1,\dots,t_m$.
If $\phi \in \cU$ satisfies
  \begin{equation}
    \label{eq:perm-inv}
 \phi(\textbf{p}(t))= \phi(t) \qquad \text{for all $t \in \R^{m}$, $\textbf{p} \in \cP$,}
  \end{equation}
 then
  \begin{equation}
    \label{eq:perm-u-inv}
 u_\phi(z,\textbf{p}(t))= u_\phi(z,t)  \qquad \text{for all $(z,t) \in \Omega$, $\textbf{p} \in \cP$}
  \end{equation}
and
  \begin{equation}
    \label{eq:perm-u-normal-inv}
\partial_{\nu_\phi} u_\phi(z,\textbf{p}(t))= \partial_{\nu_\phi} u_\phi(z,t)  \qquad \text{for all $(z,t) \in \partial \Omega$, $\textbf{p}\in \cP$.}
  \end{equation}
\end{itemize}
\end{Lemma}

\proof
Let
$$
X:= \{u \in C^{2,\alpha}_{p,e}(\overline \O)\::\: u = 0 \;\text{on $\partial \O$}\}\qquad \text{and}\qquad Y= C^{0,\alpha}_{p,e}(\overline \O).
$$
Moreover, let $\cL(X,Y)$ denote the space of bounded linear operators $X \to Y$, and let $\cI(X,Y) \subset \cL(X,Y)$ denote the subset of
topological isomorphisms $X \to Y$.  Since the metric coefficients of $g_\phi$ are smooth functions of $\phi$ and $\nabla \phi$, it is easy to see that the map
$$
\Upsilon: \cU \to \cL(X,Y), \qquad \phi \mapsto \Upsilon(\phi):= -\Delta_{g_\phi}
$$
is smooth. Moreover, for $\phi \in \cU$, the definition of $g_\phi$
implies that $\Delta_{g_\phi}$ is an elliptic, coercive second order
differential operator in divergence form with
$C^{1,\alpha}(\overline \O)$-coefficients. This immediately implies
that, by the maximum principle and elliptic regularity,  $\Upsilon(\phi) \in
\cI(X,Y)$ for every $\phi \in \cU$, and consequently the problem
(\ref{eq:Proe1-1}) has a unique solution $u_\phi \in X$ for every
$\phi \in \cU$. We now recall that $\cI(X,Y) \subset \cL(X,Y)$ is an
open set and that the inversion
$$
{\rm inv}: \cI(X,Y) \to \cI(Y,X), \qquad {\rm inv}(A)= A^{-1}
$$
is smooth. Since $u_\phi= \inv(\Upsilon(\phi)) 1$, the smoothness of the map in (\ref{eq:def-smooth-map}) follows.\\
Next, to show (i), we fix $\phi \in \cU$ and note that $u_\phi= \tilde u_\phi \circ \Psi_\phi$, where $\tilde u_\phi$ is the unique solution of the problem
$$
-\Delta \tilde u_\phi = 1 \quad \text{in $\O_\phi$,}\qquad \tilde u_\phi= 0\quad \text{on $\partial \O_\phi$.}
$$
Since $\O_\phi$ is invariant under rotations in the $z$-variable, the
uniqueness implies that $\tilde u_\phi$ is radially symmetric in $z$
and hence $u_\phi$ is also radially symmetric by the definition of
$\Psi_\phi$.  Moreover, the outer unit normal $\mu_\phi:
\partial \O_\phi \to \R^{n+m}$ with respect to $g_{eucl}$ is
equivariant with respect to rotations in $z$ by (\ref{eq:def-mu-phi}), i.e.,
$$
\mu_\phi(Az,t)= \tilde A \mu(z,t) \qquad \text{for all $(z,t) \in \partial \O_\phi$, $A \in O(n)$,}
$$
where $\tilde A \in O(n+m)$ is defined by $\tilde A(z,t)= (Az,t)$. It then follows that the function $\partial_{\mu_\phi} \tilde u_\phi: \partial \O_\phi \to \R$ is also radially symmetric in $z$, whereas by (\ref{eq:rel-mu-phi-nu-phi}) we have
$$
\partial_{\nu_{\phi}} u_\phi =  d u_{\phi} \nu_\phi = [d \tilde u_{\phi} \circ \Psi_{\phi}]  [\mu_{\phi} \circ \Psi_{\phi}]= [\partial_{\mu_\phi} \tilde u_\phi] \circ \Psi_\phi\quad \text{on $\partial \O$.}
$$
As a consequence, the function $\partial_{\nu_{\phi}} u_\phi$ is also radially symmetric in the $z$-variable, as claimed.\\
Next we note that (ii) follows from the fact that in case $\phi \equiv  \lambda>0$ we have
\begin{equation}
  \label{eq:metric-lambda}
g_{\lambda}^{ij} = \lambda^{-2} \delta_{ij} \; \text{if $1 \le i,j
\le n$ and}\; g_\lambda^{ij}=g_\lambda^{ji}= \delta_{ij}\; \text{for
$1 \le i \le n+m$, $n+1 \le j \le n+m$}
\end{equation}
and therefore
\begin{equation}
  \label{eq:laplace-lambda}
\Delta_{g_\lambda} =  \lambda^{-2} \Delta_z + \Delta_t.
\end{equation}
Consequently, the function $(z,t) \mapsto \frac{\lambda^2-|\lambda z|^2}{2n}$ solves (\ref{eq:Proe1-1}) for $\phi \equiv \lambda$, and thus it coincides with $u_\lambda$, as claimed.\\
Finally, to show (iii), we fix $\phi \in \cU$ such that~(\ref{eq:perm-inv}) holds, and we let $\tilde u_\phi$ be defined as in the proof of (i). Then $\O_\phi$ is invariant under coordinate permutations, and by uniqueness this implies that
$$
 \tilde u_\phi(z,\textbf{p}(t))= \tilde u_\phi(z,t)  \qquad \text{for all $(z,t) \in \Omega_\phi$, $\textbf{p} \in \cP$.}
$$
By definition of $\Psi_\phi$, we then conclude that the function $u_\phi$ satisfies (\ref{eq:perm-u-inv}), as claimed. Moreover, by a similar argument as in the proof of (i), we find that (\ref{eq:perm-u-normal-inv}) holds as well.
\QED

By Lemma~\ref{sec:peri-solut-serr-1}, condition (\ref{eq:Proe1-1-neumann}) is equivalent to
\begin{equation}
  \label{eq:reformulated-overd}
[\partial_{\nu_{\phi}} u_{\phi}](e_1,t)=-c \qquad \text{for $t \in \R^m$,}
\end{equation}
where $e_1 \in \R^n$ is the first coordinate vector. It follows from (\ref{eq:rel-mu-phi-nu-phi}) and (\ref{eq:def-mu-phi}) that the map
$$
\cU  \to C^{1,\alpha}_{p,e}(\partial \O,\R^{n+m}), \qquad \phi \mapsto \nu_\phi
$$
is smooth, and thus we have a smooth map
\begin{equation}
  \label{eq:def-H}
H: \cU \to C^{1,\alpha}_{p,e}(\R^m), \qquad H(\phi)(t)= \partial_{\nu_{\phi}} u_{\phi} (e_1,t),
\end{equation}
whereas (\ref{eq:reformulated-overd}) writes as
\begin{equation}
  \label{eq:reformulated-overd-H}
H(\phi)\equiv -c  \qquad \text{on $\R^m$}.
\end{equation}
In order to find solutions of the latter equation bifurcating from the trivial branch of solutions $\phi \equiv \lambda$, $\lambda>0$, we need to study the linearization of $H$ at constant functions.
The following is the main result of this section.

\begin{Proposition}
\label{sec:peri-solut-serr-2}
At a constant function $\lambda >0$, the operator $\cH_\lambda : = D H(\lambda) \in \cL(C^{2,\alpha}_{p,e}(\R^m), C^{1,\alpha}_{p,e}(\R^m))$ is given by
\begin{equation}
  \label{eq:diff-H-express}
[\cH_\lambda \omega] (t) = \frac{1}{n}\Bigl( \partial_{\nu}
\psi_{\omega,\lambda}(e_1,t)-\omega(t)\Bigr) \qquad \text{for
$\omega \in C^{2,\alpha}_{p,e}(\R^m)$, $t \in \R^m$,}
\end{equation}
where $\psi_{\omega,\lambda} \in C^{2,\alpha}_{p,e}(\overline \O)$ is the unique solution of the problem
\begin{equation}
  \label{eq:psi-omega}
\left \{
  \begin{aligned}
 \Delta_z \psi_{\omega,\lambda}(z,t) &+ \lambda^2 \Delta_t \psi_{\omega,\lambda}(z,t) =0 && \qquad (z,t) \in \O,\\
\psi_{\omega,\lambda}(z,t)&= \omega(t)  && \qquad (z,t) \in \partial \O
  \end{aligned}
\right.
\end{equation}
and $\nu$ is the outer unit normal on $\partial \O$ with respect to $g_{eucl}$ given by $\nu(z,t)= (z,0)$.
\end{Proposition}

The remainder of this section is devoted to the proof of Proposition~\ref{sec:peri-solut-serr-2}. In the following, we put
$$
\tilde \nu_\phi(\omega):= [D_\phi \nu_\phi]\omega \in C^{1,\alpha}_{p,e}(\partial \O,\R^{n+m}) \qquad \text{for $\phi \in \cU,\: \omega \in C^{2,\alpha}_{p,e}(\R^m)$.}
$$
We start with the following simple observations.

\begin{Lemma}$ $\\
\label{sec:peri-solut-serr-3}
\begin{itemize}
\item[(i)] The map $\cU \to C^{2,\alpha}_{p,e}(\overline \O,\R^{n+m})$, $\phi \mapsto \Psi_\phi$ is smooth (recall \eqref{eq:def-diffeom-Psi_phi}). Moreover, the derivative $D_\phi \Psi_{\phi} \in \cL \bigl(C^{2,\alpha}_{p,e}(\R^m), C^{2,\alpha}_{p,e}(\overline \O,\R^{n+m})\bigr)$ is given by
\begin{equation}
  \label{eq:diffPsi}
[D_\phi \Psi_{\phi}]\omega(z,t) = (\omega(t)z,0) \qquad \text{for $(z,t) \in \overline \O$, $\omega \in C^{2,\alpha}_{p,e}(\R^m)$.}
\end{equation}
\item[(ii)]  Let
$$
g: \cU \to C^{2,\alpha}(\overline \O), \qquad \phi \mapsto g_\phi
$$
be a smooth map. Then the map
$$
\cG: \cU \to C^{1,\alpha}(\partial \O), \qquad \cG(\phi)=
\partial_{\nu_{\phi}} g_\phi
$$
is smooth as well and satisfies
$$
D_\phi \cG(\phi)\omega = \partial_{\tilde \nu_\phi (\omega)} g_\phi +  \partial_{\nu_{\phi}} \Bigl([D_{\phi} g_\phi]\omega\Bigr) \qquad \text{for $\omega \in C^{2,\alpha}_{p,e}(\R^m)$.}
$$

\end{itemize}
 \end{Lemma}

\proof
(i) follows immediately from the definition of $\Psi_\phi$.\\
(ii) For $\phi \in \cU$ and $(z,t) \in \partial \Omega$ we have
$$
\cG(\phi)(z,t) = \partial_{\nu_\phi} g_\phi (z,t)= d g_\phi(z,t) \nu_{\phi}(z,t).
$$
Hence $\cG$ is smooth as a bilinear form composed with two smooth functions, and its
derivative is given by
$$
D_\phi \cG(\phi)\omega = d  [(D_\phi g_\phi)\omega]  \nu_{\phi} + d g_\phi  \bigl(\tilde \nu_\phi (\omega)\bigr) = \partial_{\nu_{\phi}} \Bigl([D_{\phi} g_\phi]\omega\Bigr)
+ \partial_{\tilde \nu_\phi (\omega)} g_\phi,
$$
as claimed.
\QED

Next we consider the function
$$
\bar u: \R^n \times \R^m \to \R,\qquad \bar u(z,t)= -\frac{|z|^2}{2n},
$$
and the smooth map
$$
\cU  \to C^{2,\alpha}(\overline \O),\qquad \phi \mapsto u^\phi:= \bar u \circ \Psi_\phi.
$$
Since $\bar u$ satisfies $-\Delta \bar u = 1$ in $\R^n \times \R^m$, for every $\phi \in \cU$ we have
\begin{equation}
\label{eq:Proe1-1-0}
    -\Delta_{g_\phi} u^\phi= 1 \qquad \textrm{in $\O$}.
\end{equation}

\begin{Lemma}
 \label{sec:peri-solut-serr}
The map
$$
\cU  \to C^{1,\alpha}(\R^m), \quad \phi \mapsto
\partial_{\nu_\phi} u^\phi(e_1,\cdot)
$$
is smooth, and its derivative at a constant function $\phi \equiv  \lambda>0$ satisfies
\begin{equation}
  \label{eq:deriv-form}
  \bigl[D_\phi |_{\phi= \lambda}\, \partial_{\nu_\phi} u^\phi(e_1,\cdot) \bigr]\omega(t)=-\frac{1}{n} \omega(t)  \qquad \text{for $\omega \in C^{2,\alpha}_{p,e}(\R^m)$ and $t \in \R^m$.}
\end{equation}
 \end{Lemma}

 \proof
The smoothness follows directly from Lemma~\ref{sec:peri-solut-serr-3}(ii). To see (\ref{eq:deriv-form}), we consider the smooth function
$$
\cU \to C^{1,\alpha}_{p,e}(\R^{m}),\qquad \phi
\mapsto\tilde{\mu}_\phi(\cdot) := \mu_\phi(\Psi_\phi (e_1,\cdot)),
$$
where $\mu_\phi: \partial \O_{\phi} \to \R^{n+m}$ is the unit outer normal with respect to $g_{eucl}$ given in (\ref{eq:def-mu-phi}). Let $\omega \in C^{2,\alpha}_{p,e}(\R^m)$. Since the map $\cU \to C^{1,\alpha}_{p,e}(\R^{m}), \;\phi \mapsto |\tilde{\mu}_\phi(\cdot)|^2_{g_{eucl}}$ is constant, we have
\begin{equation}
  \label{eq:orthogonal-mu-vp}
0 = [D_\phi\, |\tilde{\mu}_\phi(\cdot)|^2_{g_{eucl}}] \omega= 2
\langle [D_\phi \, \tilde{\mu}_\phi(\cdot)]\omega
,\tilde{\mu}_\phi(\cdot)\rangle_{g_{eucl}}\quad\textrm {on}\quad \R^m.
\end{equation}
Moreover, by (\ref{eq:rel-mu-phi-nu-phi}) we have
$\tilde{\mu}_\phi(t)=d\Psi_\phi(e_1,t)\nu_\phi(e_1,t)$ for $t\in
\R^m$ and therefore
\begin{align*}
\partial_{\nu_\phi} u^{\phi}(e_1,\cdot) = d u^\phi(e_1,\cdot) \nu_{\phi} (e_1,\cdot) =d \bar u (\Psi_{\phi}(e_1,\cdot))  &d \Psi_{\phi}(e_1,\cdot)  \nu_{\phi}(e_1,\cdot)= d \bar u (\Psi_{\phi}(e_1,\cdot))  \tilde{\mu}_\phi(\cdot)\\
  &= \langle \tilde{\mu}_\phi(\cdot) , \nabla_{g_{eucl}} \bar u(\Psi_{\phi}(e_1,\cdot)) \rangle_{g_{eucl}}
\end{align*}
on $\R^m$. Consequently,
\begin{align}
[D_\phi |_{\phi= \lambda} \, \partial_{\nu_\phi} u^{\phi}
(e_1,\cdot)]\omega &= \langle  [D_\phi |_{\phi= \lambda} \,
\tilde\mu_\phi (\cdot)]\omega , \nabla_{g_{eucl}} \bar u (\Psi_{\lambda}(e_1,\cdot)) \rangle_{g_{eucl}}\\
& +
\langle  \tilde{\mu}_\lambda(\cdot)  , [D_\phi |_{\phi= \lambda}\,  [ \nabla_{g_{eucl}} \bar u (\Psi_{\phi}(e_1,\cdot))]\omega  \rangle_{g_{eucl}} \nonumber\\
&=: I_1(\omega)+ I_2(\omega). \label{I-1-I-2}
\end{align}
Since, by definition, $\nabla_{g_{eucl}} \bar u(z,t)=- \frac{1}{n}(z,0)$ for
$(z,t) \in \R^n \times \R^m$, we get
\begin{equation}
  \label{eq:formula-nabla-bar-u-circ}
\nabla_{g_{eucl}} \bar u (\Psi_{\phi}(e_1,t)) =
 - \frac{1}{n}(\phi(t) e_1,0)\qquad \text{for $t \in \R^m$}
\end{equation}
and thus
$$
[D_\phi |_{\phi= \lambda}   [(\nabla_{g_{eucl}}\bar u(\Psi_{\phi}(e_1, \cdot))]\omega(t)= -\frac{1}{n} (\omega(t)e_1,0)
\qquad t \in \R^m.
$$
Noting also that $\tilde {\mu}_\lambda \equiv (e_1,0) \in \R^{n+m}$ on $\R^m$, we infer that
\begin{equation}
  \label{eq:I2}
I_2(\omega)(t) = \langle  \tilde{\mu}_\lambda(t)  , -\frac{1}{n}
(\omega(t)e_1,0) \rangle_{g_{eucl}} = -\frac{1}{n}\omega(t) \qquad
\text{for $t \in \R^m$}.
\end{equation}
From (\ref{eq:formula-nabla-bar-u-circ}) we also deduce that
$$
(\nabla_{g_{eucl}}\bar u (\Psi_{\lambda}(e_1,\cdot)) \equiv-
\frac{1}{n}(\lambda e_1,0)\equiv - \frac{\lambda}{n}
\tilde{\mu}_\lambda(\cdot),
$$
so the identity (\ref{eq:orthogonal-mu-vp}) with $\phi = \lambda$
implies that
\begin{equation}
  \label{eq:I1}
I_1(\omega)=0 \qquad \text{on $\R^m$.}
\end{equation}
Now (\ref{eq:deriv-form}) follows by combining (\ref{I-1-I-2}), (\ref{eq:I2}) and (\ref{eq:I1}).
 \QED

We are now in a position to complete the\\

{\bf Proof of Proposition~\ref{sec:peri-solut-serr-2}.}
For $\phi \in \cU$, we note that the function $a_\phi:= u^\phi-u_\phi \in C^{2,\alpha}_{p,e}(\overline \Omega)$ satisfies
\begin{equation}
\label{eq:Proe1-2}
 \left \{
 \begin{aligned}
    -\Delta_{g_\phi} a_\phi&=0 &&\qquad \textrm{ in $\O$}\\
a_\phi&= u^\phi \qquad &&\qquad \textrm{on $\partial \O$}.
  \end{aligned}
\right.
\end{equation}
Moreover,  in the case where $\phi \equiv \lambda>0$, we have $u^\lambda(z,t)= -\frac{|\lambda z|^2}{2n}$ and therefore
\begin{equation}
  \label{eq:special-case-lambda}
a_\lambda(z,t) = -\frac{\lambda^2}{2n}\qquad \text{for $(z,t) \in \overline \Omega$}
 \end{equation}
by Lemma~\ref{sec:peri-solut-serr-1}(ii). Now, consider the smooth map $T: \cU \to C^{0,\alpha}(\overline \Omega)$ given by
$$
 T(\phi) = \Delta_{g_\phi} a_\phi = g_{\phi}^{ij} \partial_{ij} a_{\phi} + \frac{1}{\sqrt{|g_\phi|}} \partial_i \Bigl(\sqrt{|g_\phi|}g_{\phi}^{ij}\Bigr)\partial_j
 a_\phi.
$$
By (\ref{eq:Proe1-2}) we have $T \equiv 0$ on $\cU$. Thus for every $\omega \in C^{2,\alpha}_{p,e}(\R^m)$ we have
\begin{equation}
  \label{eq:0-identity}
0= D T(\phi) \omega =  \Delta_{g_\phi} [D_\phi a_\phi]\omega   + h_{\phi}^{ij} \partial_{ij} a_{\phi} + \ell^j_\phi \partial_j a_\phi
\end{equation}
with
$$
h_{\phi}^{ij}:= [D_\phi  g_{\phi}^{ij}]\omega  \qquad \text{and}\qquad \ell^j_\phi := \bigl[D_\phi \frac{1}{\sqrt{|g_\phi|}} \partial_i \bigl(\sqrt{|g_\phi|}g_{\phi}^{ij}\bigr)\bigr] \omega.
$$
Evaluating (\ref{eq:0-identity}) at $\phi= \lambda$ and using that the function $a_\lambda$ is constant in $\Omega$ by (\ref{eq:special-case-lambda}), we find that the function
$$
\tau_{\omega,\lambda}:= [D_\phi \big|_{\phi=\lambda} a_\phi] \omega \in C^{2,\alpha}_{p,e}(\Omega)
$$
satisfies
$$
 \Delta_z \tau_{\omega,\lambda}+ \lambda^2 \Delta_t \tau_{\omega,\lambda} = \lambda^2 \Delta_{g_\lambda} \tau_{\omega,\lambda} =0 \qquad \text{in $C^{0,\alpha}(\overline \Omega)$.}
$$
(here the first equality follows from (\ref{eq:laplace-lambda})). Moreover, differentiating the boundary condition in (\ref{eq:Proe1-2}) and using Lemma~\ref{sec:peri-solut-serr-3}(i) gives
\begin{align*}
\tau_{\omega,\lambda}(z,t) &= [D_\phi \big|_{\phi=\lambda}\, \bar u
\circ \Psi_\phi]\omega(z,t) = d \bar u (\Psi_\lambda(z,t))
[D_\phi \big|_{\phi=\lambda} \Psi_\phi]\omega(z,t) \\
&= d \bar u (\lambda z, t)(\omega(t)z,0)
=- \frac{\lambda}{n}\omega(t) \qquad \text{for $(z,t)
\in \partial \O$.}
\end{align*}
By Lemma~\ref{sec:peri-solut-serr-3}(ii) and since $a_\lambda$ is constant in $\overline \Omega$, we also have that
\begin{equation}
  \label{eq:proof-linearization-1}
[D_\phi \big|_{\phi=\lambda}   \partial_{\nu_{\phi}} a_\phi]\omega  = \partial_{\tilde \nu_{\phi}(\omega)} a_\lambda + \partial_{\nu_\lambda} \tau_{\omega,\lambda} = \frac{1}{\lambda} \partial_\nu \tau_{\omega,\lambda} \qquad \text{on $\partial \Omega$},
\end{equation}
where $\nu$ is the outer unit normal on $\partial \O$ with respect
to $g_{eucl}$ given by $\nu(z,t)= (z,0)$. Combining
(\ref{eq:proof-linearization-1}) with
Lemma~\ref{sec:peri-solut-serr}, we thus find that
$$
 [D_\phi \big|_{\phi=\lambda}    H(\phi)]\omega(t)=  [D_\phi \big|_{\phi=\lambda}   \partial_{\nu_{\phi}} (u^\phi - a_\phi)(e_1,\cdot)]\omega(t)
= -  \frac{1}{n} \omega(t)-\frac{1}{\lambda} \partial_{\nu} \tau_{\omega,\lambda}(e_1,t)
$$
for $t \in \R^m$. Putting $\psi_{\omega,\lambda}:= -\frac{n}{\lambda}\tau_{\omega,\lambda}$, we then see that (\ref{eq:diff-H-express}) and (\ref{eq:psi-omega}) hold, as claimed.
\QED

\section{Spectral properties of the linearization}
\label{sec:spectr-prop-line-1}
In this section we study the spectral properties of the linearized operators  $\cH_\lambda = D H(\lambda) \in \cL(C^{2,\alpha}_{p,e}(\R^m), C^{1,\alpha}_{p,e}(\R^m))$, $\lambda>0$ considered in Proposition~\ref{sec:peri-solut-serr-2}. We start with the following observation.

\begin{Proposition}\label{proposition-eigenvalues}
Let $\lambda>0$. The functions
\begin{equation}
  \label{eq:def-omega_k}
\omega_{k} \in C^{2,\alpha}_{p,e}(\R^m), \quad \omega_{k}(s)=\prod_{j=1}^m \cos(k_j s_j), \quad k \in \bigl(\N \cup \{0\}\bigr)^m
\end{equation}
are eigenfunctions of $\cH_\lambda$ in the sense that
\begin{equation}
  \label{eq:eigenv-eq}
[\cH_\lambda \omega_k](t)= \sigma(|k| \lambda)\, \omega_k(t) \qquad
\text{for $t \in \R^{m}$ with $ |k|=
\sqrt{k_1^2+ \dots + k_m^2}$.}
\end{equation}
Here the function $\sigma$ is defined
by \be \label{eq:lkhk} \sigma: [0,\infty) \to \R,\qquad
\sigma(\rho)=\frac{1}{n}\Bigl(\rho \frac{h'(\rho)}{h(\rho)} -
1\Bigr),
 \ee
where $h : [0,\infty) \to \R$ is the unique solution of the initial
value problem
\begin{equation}\label{eq:problemr-0}
\begin{cases}
\displaystyle h''(\rho)+\frac{n
-1}{\rho}h'(\rho)-h(\rho)=0 \\
\displaystyle h(0)=1, \quad h'(0)=0.
\end{cases}
\end{equation}
Furthermore, in case $n=1$ we have
\begin{equation}
  \label{eq:sigma-k-N=1}
 \sigma(\rho) = \rho \tanh (\rho)-1,
\end{equation}
and in case $ n \ge 2$ we have
\begin{equation}
  \label{eq:sigma-k-N2}
 \sigma(\rho) = \frac{1}{n}\Bigl(\rho \frac{I_{\nu+1}(\rho)}{I_\nu(\rho) }-1 \Bigr),
\end{equation}
where $I_\nu$ is the modified Bessel function of the first kind of order $\nu= \frac{n-2}{2}$.
\end{Proposition}

\proof Fix $\lambda >0$ and $k \in \bigl(\N \cup \{0\}\bigr)^m$. For
$\omega = \omega_k \in C^{2,\alpha}_{p,e}(\R^m)$ as defined in
(\ref{eq:def-omega_k}), the unique solution
$\psi_{\omega_k,\lambda}$ of (\ref{eq:psi-omega}) is given by
$\psi_{\omega_k, \l}(z,t)=b(|z|) \omega_k(t)$, where the function
$b:[0,1] \to \R$ is the unique solution to
\begin{equation}\label{eq:problem4}
\begin{cases}
\displaystyle b''+\frac{n
-1}{r}b'- \l^2 |k|^2b=0,\qquad r\in(0,1) \\
\displaystyle \displaystyle b'(0)=0,\quad  b(1)=1.
\end{cases}
\end{equation}
Hence we have that
$$
[\cH_\lambda \omega_k](e_1,t)= \frac{1}{n}\Bigl( \partial_{\nu}
\psi_{\omega_k,\lambda}(e_1,t)-\omega_k(t)\Bigr)=
\frac{1}{n}(b'(1)-1)\omega_k(t) \qquad \text{for $t \in \R^m$}.
$$
Now putting $\rho_0 = \lambda |k|$ and considering $\tilde h: [0,\rho_0] \to \R$ defined by $\tilde h(\rho):=b(\frac{\rho}{\rho_0})$, we see that $\tilde h$ satisfies
\begin{equation}\label{eq:problemr}
\begin{cases}
\displaystyle \tilde h''(\rho)+\frac{n
-1}{\rho}\tilde h'(\rho)-\tilde h(\rho)=0,\qquad \rho\in(0,\rho_0),\\
\displaystyle \displaystyle \tilde h'(0)=0,\qquad \tilde h(\rho_0)=1.
\end{cases}
\end{equation}
Consequently, $\tilde h = \frac{h}{h(\rho_0)}$ in $[0,\rho_0]$, where
$h : [0,\infty) \to \R$ is the unique solution of the initial value problem (\ref{eq:problemr-0}). Moreover,
$$
\frac{1}{n}(b'(1)-1) = \frac{1}{n} \bigl(\rho_0 \tilde h'(\rho_0) -1 \bigr)=  \frac{1}{n}\bigl(\rho_0 \frac{h'(\rho_0)}{h(\rho_0)} - 1\bigr) = \sigma(\rho_0),
$$
as claimed in \eqref{eq:lkhk}. Now in case $n=1$ we have  $h(\rho)= \cosh(\rho)$ for $\rho>0$ and thus \eqref{eq:lkhk} follows.\\
In case $n \ge 2$, we consider $g(\rho):=\rho^{\nu}h(\rho)$ with $\nu:=\frac{n-2}{2}$, so that (\ref{eq:problemr-0}) transforms into
the following (modified) Bessel
equation:
\begin{equation*}
\displaystyle
g''(\rho)+\frac{1}{\rho}g'(\rho)-\biggl(1+\frac{\nu^{2}}{\rho^{2}}\biggl)g(\rho)=0.
\end{equation*}
 Up to a constant, the unique locally bounded solution to this equation is the modified Bessel function of the first kind $I_\nu$. Since $I_\nu>0$ on $(0,\infty)$, we thus have
 \begin{equation}
   \label{h-I-nu-relation}
 h(\rho)= c\, \rho^{-\nu}I_\nu(\rho) \qquad \text{for $\rho>0$ with a constant $c>0$.}
 \end{equation}
In fact, it will follow from \eqref{eq:problemr-0} and (\ref{Bges}) below that $c= 2^{\nu}\Gamma(\nu+1)$, but we do not need this. Using (\ref{h-I-nu-relation}) together with the recurrence formula $\rho I'_{\nu}(\rho)-\nu
 I_{\nu}(\rho)=\rho I_{\nu+1}(\rho)$ (see e.g.  \cite[Section 7.11]{El}), we find that
$$
\frac{h'(\rho)}{h(\rho)}= \frac{I_\nu'(\rho)-\frac{\nu}{\rho} I_\nu(\rho)}{I_\nu(\rho)} = \frac{I_{\nu+1}(\rho)}{I_{\nu}(\rho)}\qquad \text{for $\rho>0$.}
$$
Therefore \eqref{eq:lkhk} yields
$$
\s(\rho)=\frac{1}{n}\Bigl(\rho \frac{I_{\nu+1}(\rho)}{I_{\nu}(\rho)} - 1\Bigr),
$$
 as claimed.
\QED

 The following lemma
gives the asymptotic behavior of the function $\sigma$.
 \begin{Lemma}\label{eq:asym}
The function $\sigma:(0,\infty)\longrightarrow \mathbb{R}$ has the
following asymptotic properties:
\begin{itemize}
\item[(i)]  $\lim \limits_{\rho \to \infty} \frac{\sigma(\rho)}{\rho} = \dfrac{1}{n},$
\item[(ii)]  $\lim \limits_{\rho \to 0} \sigma(\rho) = -\dfrac{1}{n}.$
\end{itemize}
 \end{Lemma}

 \proof In case $n=1$, both (i) and (ii) follow immediately from (\ref{eq:sigma-k-N=1}). In case $n \ge 2$, (i) follows from \eqref{eq:sigma-k-N2} and the asymptotic formula
$$
\lim_{\rho\rightarrow +\infty}\frac{I_{\tau}(\rho)}{\frac{1}{\sqrt{2\pi\rho }}e^{\rho}}=1\qquad \qquad \text{for every $\tau \ge 0$,}
$$
see \cite[Section 7.13.1]{El}. Moreover  (see e.g. \cite[Section 7.2.2]{El}), for $\tau \ge 0$, we have the power series  expression of $I_{\tau}$ given by
 \begin{equation}\label{Bes}
 I_{\tau}(\rho)=\sum^{\infty}_{i=0}\frac{(\frac{1}{2}\rho)^{\tau+2i}}{i!\Gamma(\tau+i+1)}, \qquad \textrm{ for $\t\geq 0$ and $\rho>0$}.
 \end{equation}
We can write
\begin{equation}\label{Bges}
 I_{\tau}(\rho)=(\frac{1}{2}\rho)^{\tau}\biggl[\frac{1}{\Gamma(\tau+1)}+\sum^{\infty}_{i=1}\frac{(\frac{1}{2}\rho)^{2i}}{i!\Gamma(\tau+i+1)}\biggl],
 \end{equation}
which shows that $\lim \limits_{\rho \to 0} \rho\frac{I_{\tau+1}(\rho)}{I_{\tau}(\rho)}=0$ for every $\tau \ge 0$. Together with (\ref{eq:sigma-k-N2}) this gives (ii).\QED

Next, we show that the functions $\lambda \mapsto \sigma(\lambda)$ are
strictly increasing on $(0,\infty)$.

\begin{Lemma}\label{eq: Deri}
We have $\sigma'(\rho)>0$ for $\rho>0$. Moreover, $\sigma$ has
exactly one zero in $(0,\infty)$.
\end{Lemma}
\proof By Lemma~\ref{eq:asym}, we only need to show that
$\sigma'(\rho)>0$ for $\rho>0$. In case $n=1$,
(\ref{eq:sigma-k-N=1}) gives that $\sigma'(\rho)=
\frac{\rho}{\cosh^2(\rho)} + \tanh(\rho) >0$ for $\rho>0$. In case
$n \ge 2$, we use \eqref{eq:lkhk} and calculate that
$$
n h^2(\rho) \sigma'(\rho) = h'(\rho) h(\rho) + \rho \Bigl(h''(\rho)h(\rho)- h'(\rho)^2\Bigr)
= (2-n) h'(\rho)h(\rho) + \rho (h^2(\rho)- h'(\rho)^2)
$$
for $\rho > 0$. For the latter equality, we used the fact that
\begin{equation}
  \label{eq:alternative-proof-mono-1}
h''(\rho) = h(\rho) +\frac{1-n}{\rho} h'(\rho) \qquad \text{for $\rho>0$}
\end{equation}
as a consequence of (\ref{eq:problemr-0}). It then suffices to show that the function
$$
\rho \mapsto  j(\rho):= \rho^{n-1} n h^2(\rho) \sigma'(\rho)= (2-n)\rho^{n-1}h'(\rho)h(\rho)
+ \rho^n (h^2(\rho)- h'(\rho)^2)
$$
is positive on $(0,\infty)$. Since $j(0)=0$, it suffices to show that $j'(\rho)>0$ for $\rho>0$. Using (\ref{eq:alternative-proof-mono-1}) again, we find that
\begin{align*}
j'(\rho)&=  (n-1) (2-n)\rho^{n-2}h'(\rho)h(\rho) +(2-n)\rho^{n-1} \bigl(h''(\rho)h(\rho)+h'(\rho)^2\bigr)\\
&+ n \rho^{n-1} \bigl(h^2(\rho)- h'(\rho)^2\bigr) +2 \rho^{n}\bigl(h(\rho)h'(\rho) -h'(\rho)h''(\rho)\bigr)\\
&=  (n-1) (2-n)\rho^{n-2}h'(\rho)h(\rho) +(2-n)\rho^{n-1} \bigl(h^2(\rho)+\frac{1-n}{\rho} h'(\rho)h(\rho)+h'(\rho)^2\bigr)\\
&+ n \rho^{n-1} \bigl(h^2(\rho)- h'(\rho)^2\bigr) +2(n-1) \rho^{n-1}h'(\rho)^2\\
&= \rho^{n-1} \Bigl(2 h^2(\rho)  + (2-n) h'(\rho)^2  -n h'(\rho)^2+ 2(n-1)h'(\rho)^2\Bigr)\\
&= 2\rho^{n-1} h^2(\rho)= 2c^2 \rho I_{\nu}^2(\rho)  \qquad \text{for $\rho>0$}
\end{align*}
with $\nu= \frac{n-2}{2}$ and $c>0$ as in (\ref{h-I-nu-relation}). Since $I_\nu(\rho)>0$ for $\rho>0$, we thus conclude that $j'(\rho)>0$ for $\rho>0$, as required.
\QED

In the following, we consider the Sobolev spaces
\begin{equation}
  \label{eq:def-hpe}
H^{j}_{p,e} := \Bigl \{v \in H^{j}_{loc}(\R^m) \::\: \text{$v$ even,
$2\pi$-periodic in $t_1,\dots,t_m$}\Bigl \}, \qquad j \in \N \cup
\{0\},
\end{equation}
and we put $L^2_{p,e}:= H^{0}_{p,e}$. Note that $L^2_{p,e}$ is a Hilbert space with scalar product
$$
(u,v) \mapsto \langle u,v \rangle_{L^2} := \int_{[0,2\pi]^m} u(t)v(t)\,dt \qquad \text{for $u,v \in L^2_{p,e}$.}
$$
We denote the induced norm by $\|\cdot\|_{L^2}$. For the functions $\omega_{k}$ in (\ref{eq:def-omega_k}) we then have
$$
\|\omega_0\|_{L^2}=(2\pi)^{m/2}, \qquad \quad  \|\omega_k\|_{L^2}=\pi^{m/2} \quad \text{for $k \in (\N \cup \{0\})^m$, $k \not = 0$,}
$$
and the functions $\frac{\omega_{k}}{\|\omega_k\|_{L^2}}$ form an orthonormal basis for $L^2_{p,e}$. Moreover,
$H^j_{p,e} \subset L^2_{p,e}$ is characterized as the subspace of all functions $v \in L^2_{p,e}$ such that
$$
\sum_{k \in (\N \cup \{0\})^m} (1+|k|^2)^{j} \langle v,\omega_k \rangle_{L^2} ^2 < \infty.
$$
Thus, $H^j_{p,e}$ is also a Hilbert space with scalar product
\begin{equation}
  \label{eq:scp-hj}
(u,v) \mapsto \langle u,v \rangle_{H^j} :=  \sum_{k \in (\N \cup \{0\})^m} (1+|k|^2)^{j} \langle u, \omega_k \rangle_{L^2}  \langle v,\omega_k \rangle_{L^2}  \qquad \text{for $u,v \in H^j_{p,e}$.}
\end{equation}
In the following, we also consider the subspaces
\begin{equation}
  \label{eq:def-Vell}
V_\ell:= \langle \omega_k \;:\; |k|=\ell \rangle \,\subset \,\bigcap_{j \in \N} H^j_{p,e},
\end{equation}
the corresponding $\langle \cdot,\cdot \rangle_{L^2}$-orthogonal projections $P_\ell: L^2_{p,e} \to L^2_{p,e}$ on $V_\ell$, and the
complements
\begin{equation}
  \label{eq:def-Wjell}
Z^j_\ell := \{v \in H^j_{p,e}\::\: P_\ell v = 0\} \:\subset \: H^j_{p,e},\qquad \ell \in \N \cup \{0\}.
\end{equation}
Since the latter spaces are closed subspaces of $H^j_{p,e}$, they are also Hilbert spaces with respect to the scalar product in (\ref{eq:scp-hj}).

\begin{Proposition}
\label{sec:spectr-prop-line} For fixed $\lambda>0$, the linear map
$\cH_{\lambda}$ defined in \eqref{eq:diff-H-express} extends to a
continuous linear map
$$
\cH_\lambda: H^{2}_{p,e} \to H^{1}_{p,e},\qquad \cH_\lambda v =
\sum_{\ell \in \N \cup \{0\}} \sigma(\lambda \ell) P_\ell v.
$$
Moreover, for any $\ell \in \N \cup \{0\}$, the operator
$$
\cH_\lambda - \sigma(\lambda \ell) \id \;:\;  Z^2_\ell \to Z^1_\ell \qquad \text{is an isomorphism.}
$$
\end{Proposition}

\proof
This follows from Proposition~\ref{proposition-eigenvalues}, Lemma~\ref{eq:asym}, Lemma~\ref{eq: Deri} and the remarks above.
\QED

\begin{Remark}
  \label{sec:spectr-remark-1}
The extension $\cH_\lambda: H^{2}_{p,e} \to H^{1}_{p,e}$ given in Proposition~\ref{sec:spectr-prop-line} can be characterized as follows. For $k \in \N \cup \{0\}$ and $1 \le p \le \infty$, we consider the space
\begin{equation}
  \label{eq:def-*-spaces}
W^{k,p}_{b}(\Omega):= \{ \psi \in W^{k,p}_{\loc}(\Omega)\::\: \text{$\psi \in W^{k,p}(\Omega')$ for every bounded subset $\Omega' \subset \Omega$}\}
\end{equation}
Given $\omega \in H^2_{p,e}$, standard elliptic theory shows that there is a unique solution $\psi \in W^{2,2}_{b}(\Omega)$ of the problem
\begin{equation}
  \label{eq:psi-omega-0}
\left \{
  \begin{aligned}
 \Delta_z \psi (z,t) &+ \lambda^2 \Delta_t \psi(z,t) =0 && \qquad (z,t) \in \O,\\
\psi(z,t)&= \omega(t)  && \qquad (z,t) \in \partial \O
  \end{aligned}
\right.
\end{equation}
which is even and $2\pi$-periodic in $t_1,\dots,t_m$. Then
 $\cH_\lambda \omega$ is given by
$$
[\cH_\lambda \omega] (t) = \frac{1}{n}\Bigl( \partial_{\nu}
\psi(e_1,t)-\omega(t)\Bigr) \qquad \text{for a.e. $t \in \R^m$,}
$$
where $\nu$ is the outer unit normal on $\partial \O$ with respect to $g_{eucl}$ and $\partial_{\nu}
\psi$ is considered in the sense of traces. This can be easily seen by approximating $\omega$ in $H^2_{p,e}$ with functions in $C^{2,\alpha}_{p,e}(\R^m)$ and using standard elliptic estimates.
\end{Remark}

\section{Proof of Theorem~\ref{teo1}}
\label{sec:compl-proof-theor}

In the following, we let $\cP \subset \cL(\R^m)$ denote the subset of all coordinate permutations, and we consider the spaces
\begin{align*}
&X:= \{\vp\in C^{2,\alpha}_{p,e}(\R^m)\::\: \vp(t)= \vp(\textbf{p}(t))  \text{ for all $t \in \R^m$, $\textbf{p} \in \cP$}\},\\
&Y:= \{\vp \in C^{1,\alpha}_{p,e}(\R^m)\::\: \vp(t)= \vp(\textbf{p}(t))  \text{ for all $t \in \R^m$, $\textbf{p} \in \cP$}\}.
\end{align*}
We also consider the nonlinear operator $H$ defined in (\ref{eq:def-H}), and we note that $H$ maps $\cU \cap X$ into $Y$
by Lemma~\ref{sec:peri-solut-serr-1}(iii). Consider the open set
\begin{equation}
  \label{eq:def-cO}
\cO:= \{(\lambda,\vp) \in \R \times X \::\: \lambda>0, \: \vp > - \lambda \}  \subset \R \times X.
\end{equation}
The proof of Theorem~\ref{teo1} will be completed by applying the Crandall-Rabinowitz Bifurcation theorem to the smooth nonlinear operator
\be \label{eq:def--G}
G: \cO \subset \R \times X  \to Y, \qquad G(\lambda, \vp) =
H(\lambda+\vp)+\frac{\lambda}{n}.
\ee
Recalling the formula of $u_\lambda$ in Lemma
\ref{sec:peri-solut-serr-1}, we have
$$
G(\lambda,0)(t) = H(\lambda)(t)+ \frac{\lambda}{n} =
\partial_{\nu_\lambda}u_\lambda(e_1,t)+\frac{\lambda}{n}=    0
\qquad \text{for $t \in \R^m$, $\lambda>0$.}
$$
Moreover,
\begin{equation}
  \label{eq:cGH-lambda}
 D_\vp G (\lambda,0)= DH(\l)\big|_X = \cH_\lambda|_X \in \cL(X,Y).
\end{equation}
We have the following.

\begin{Proposition}\label{propPhi}
There exists a unique $\lambda_*=\l_*(n)>0$ such that $\sigma(\lambda_*)=0$, where the function $\sigma$ is defined in Proposition~\ref{proposition-eigenvalues}. Moreover, the linear operator
$$
\cH_*:=\cH_{\lambda_*}\big|_X \in \cL(X,Y )
$$
has the following properties.
\begin{itemize}
\item[(i)] The kernel $N(\cH_*)$ of $\cH_*$ is spanned by the function
  \begin{equation}
    \label{eq:def-v-0}
v_0 \in X, \qquad v_0(t)= \cos(t_1) + \dots + \cos(t_m).
  \end{equation}
\item[(ii)] The range of $\cH_*$ is given by
$$
R(\cH_*)= \Bigl \{v \in Y  \::\: \int_{[0,2\pi]^m} v(t) v_0(t)\,dt = 0 \Bigr\}.
$$
\end{itemize}
Moreover,
\begin{equation}
  \label{eq:transversality-cond}
\partial_\lambda \Bigl|_{\lambda= \lambda_*} \cH_\lambda v_0 \;\not  \in \; R(\cH_*).
\end{equation}
\end{Proposition}

\proof
By Lemma~\ref{eq: Deri}, there exists a unique $\lambda_*=\l_*(n)>0$ such that $\sigma(\lambda_*)=0$, which by
Proposition~\ref{proposition-eigenvalues} is equivalent to $\cH_{\lambda_*}v_0  = 0$. We put $\cH_*:= \cH_{\lambda_*}$ in the following. Consider the subspaces
\begin{align}
X^* &:= \Bigl \{v \in  X \::\: \int_{[0,2\pi]^m} v(t) v_0(t)\,dt = 0  \Bigr\} \subset X, \label{def-X*}\\
Y^* &:= \Bigl \{v \in Y \::\: \int_{[0,2\pi]^m} v(t) v_0(t)\,dt = 0 \Bigr\} \subset Y. \nonumber
\end{align}
To show properties (i) and (ii), it clearly suffices to prove that
\begin{equation}
  \label{eq:isomorphism}
\text{$\cH_*$ defines an isomorphism between $X^*$ and $Y^*$.}
\end{equation}
To prove (\ref{eq:isomorphism}), we need to introduce further spaces.
We recall the definition of $H^j_{p,e}$ in (\ref{eq:def-hpe}) and put
\begin{equation}
  \label{eq:def-hpes}
H^{j}_{\cP} := \Bigl \{v \in H^j_{p,e}  \::\: \text{$v(\textbf{p}(t))=v(t)$ for $t \in \R^m$ and $\textbf{p} \in \cP$}\}, \qquad j \in \N \cup \{0\},
\end{equation}
noting that $X = H^{2}_{\cP}  \cap C^{2,\alpha}(\R^m)$ and $Y = H^{1}_{\cP} \cap C^{1,\alpha}(\R^m)$. Proposition~\ref{sec:spectr-prop-line} implies that $\cH_*$ defines a continuous linear operator
\begin{equation}
  \label{eq:formular-H*-perm}
\cH_* : H^2_{\cP} \to H^1_{\cP}, \qquad \cH_* v = \sum_{\ell \in \N \cup \{0\}} \sigma(\lambda_* \ell) P_\ell v,
\end{equation}
Next we put
$$
\tilde V^j:=  H^j_{\cP} \cap V_1,\quad  \tilde Z^j:=  H^j_{\cP} \cap Z_1^j \quad \subset H^j_\cP  \qquad \text{for $j=1,2$,}
$$
where the spaces $V_1$ resp. $Z_1^j$ are defined in (\ref{eq:def-Vell}) and (\ref{eq:def-Wjell}), respectively.
  We note that $\tilde V^1$ is one-dimensional and spanned by the function $v_0$ defined in (\ref{eq:def-v-0}). Since the spaces $\tilde V^j$ and $\tilde Z^j$
  are invariant with respect to coordinate permutations $\textbf{p} \in \cP$, we deduce from Proposition~\ref{sec:spectr-prop-line} and our choice of $\lambda_*$ that
\begin{equation}
  \label{eq:Sobolev-invertible}
\text{$\cH_*$ defines an isomorphism $\tilde Z^2 \to \tilde Z^1$}.
\end{equation}
Moreover, since $X^*= \tilde Z^2 \cap X$ and $Y^* = \tilde Z^1 \cap Y$, we see that $\cH_*: X^* \to Y^*$
is well defined and injective. To establish surjectivity, let $f \in Y^*$. By (\ref{eq:Sobolev-invertible}),
there exists $\omega \in \tilde Z^2 \subset H^2_{p,e}$ such that $\cH_* \omega=f$. As noted in Remark~\ref{sec:spectr-remark-1},
we then have
\begin{equation}
   \label{eq:corr-version-add-0}
\partial_\nu  \psi (\pm e_1,t) - \omega(t)= n f(t)\quad \text{for a.e. $t \in \R^m$,}
\end{equation}
where $\psi \in W^{2,2}_{b}(\Omega)$ is the unique solution of (\ref{eq:psi-omega-0}) which is even and $2\pi$-periodic in $t_1,\dots,t_m$. We claim that
\begin{equation}
  \label{eq:corr-version-add-1}
\psi \in C^{2,\alpha}_{p,e}(\overline \O).
\end{equation}
Indeed, this follows from \cite[Theorem 6.2.3.1]{grisvard} once we have shown that
\begin{equation}
  \label{eq:corr-version-add-2}
\psi \in W^{2,p}_{b}(\O) \qquad \text{for some $p > N$,}
\end{equation}
where the space $W^{2,p}_{b}(\O)$ is defined in (\ref{eq:def-*-spaces}). To see this, we show by induction that
\begin{equation}
  \label{eq:corr-version-add-3}
\psi \in W^{2,p_k}_{b}(\O)
\end{equation}
for a sequence of numbers $p_k \in [2,\infty)$ satisfying $p_0=2$ and $p_{k+1} \ge \frac{N-1}{N-2}p_{k}$ for $k \ge 0$. We already know that (\ref{eq:corr-version-add-3}) holds for $p_0=2$. So let us assume that (\ref{eq:corr-version-add-3}) holds for some $p_k \ge 2$. We distinguish two cases.\\
If $p_k < N$, then the trace theorem implies that
$$
\psi \big|_{\partial
\O}\in  W^{1,p_{k+1}}_{loc}(\partial \O) \qquad \text{with $p_{k+1}:= (\frac{N-1}{N-p_k})p_k \ge \frac{N-1}{N-2} p_k$,}
$$
so that $\omega \in W^{1,p_{k+1}}_{loc}(\R^m)$. Since also $f \in C^{1,\alpha}(\R^m) \subset W^{1,p_{k+1}}_{loc}(\R^m)$ and, by (\ref{eq:corr-version-add-0}),
$$
\partial_\nu  \psi (z,t) +\psi(z,t)=  \partial_\nu  \psi (z,t)+\omega(t)= g(t) \quad \text{for $(z,t)  \in \partial \Omega$}
$$
with
$$
g = nf +2 \omega \in W^{1,p_{k+1}}_{loc}(\R^m),
$$
we may deduce from \cite[Theorem 2.4.2.6]{grisvard} that $\psi \in W^{2,p_{k+1}}_{b}(\O)$.\\
If $p_k \ge  N$, the trace theorem implies that $W^{1,p}_{loc}(\partial \O)$ for any $p>2$, and then we may repeat the above argument with arbitrarily chosen $p_{k+1} \ge \frac{N-1}{N-2}p_{k}$ to infer again that $\psi \in W^{2,p_{k+1}}_{b}(\O)$.\\
We thus conclude that (\ref{eq:corr-version-add-2}) holds, and hence (\ref{eq:corr-version-add-1}) follows. By passing to the trace again, we then conclude that $\omega \in  C^{2,\alpha}_{p,e}(\R^m)$. Consequently, $\omega \in C^{2,\alpha}_{p,e}(\R^m) \cap \tilde Z^2 = X^*$, and thus $\cH_*: X^* \to Y^*$ is also surjective.  Hence (\ref{eq:isomorphism}) is true.\\
It remains to prove~\eqref{eq:transversality-cond}, which follows
from Lemma~\ref{eq: Deri} and the identity
\begin{equation}\label{Dlambda}
\partial_\lambda \Bigl|_{\lambda= \lambda_*} \cH_\lambda v_0 =  \partial_\lambda \Bigl|_{\lambda= \lambda_*}\sigma(\lambda) v_0 = \sigma'(\lambda_*) v_0.
\end{equation}
\QED

\noindent {\bf Proof of Theorem~\ref{teo1} (completed).}
Recalling \eqref{eq:def--G} and \eqref{eq:def-H},  we shall  apply the Crandall-Rabinowitz Bifurcation Theorem to solve the equation
\be \label{eq:final-eq-to-solve}
G(\l,\vp)=H(\lambda+\vp)+\frac{\lambda}{n}= \partial_{\nu_{\phi}} u_{\phi} (e_1,\cdot )+\frac{\lambda}{n}= 0,
\ee
where $\phi=\l+\vp \in \cU$ and  the function  $u_{\phi}  \in C^{2,\a}(\ov{\O_{\phi}})$  is the unique solution to  the Dirichlet boundary value problem
\begin{align*}
  \begin{cases}
    -\Delta_{g_\phi} u_\phi= 1& \quad \textrm{ in}\quad\O\\
u_\phi=0&  \quad\textrm{ on} \quad \partial \O,
  \end{cases}
  \end{align*}
see  Lemma \ref{sec:peri-solut-serr-1}. Once this is done,  \eqref{eq:reformulated-overd-H} follows and thus we get \eqref{eq:reformulated-overd}, which is equivalent to \eqref{eq:Proe1-1-neumann} with $c=\frac{\lambda}{n}$. \\
 To solve equation \eqref{eq:final-eq-to-solve},   we
let $\lambda_*=\l_*(n)$ be defined as in Proposition~\ref{propPhi}, and let $X^*$ be defined as in (\ref{def-X*}).
 By Proposition~\ref{propPhi} and  the Crandall-Rabinowitz Theorem (see \cite[Theorem 1.7]{M.CR}), we then find ${\e_0}>0$ and a smooth curve
$$
(-{\e_0},{\e_0}) \to \cO, \qquad s \mapsto (\lambda_s,\varphi_s)
$$
such that
\begin{enumerate}
\item[(i)] $G(\lambda_s,\varphi_s)=0$ for $s \in (-{\e_0},{\e_0})$,
\item[(ii)] $\lambda(0)= \lambda_*$, and
\item[(iii)]  $\varphi_s = s \bigl(v_0 + \mu_s\bigr)$ for $s \in (-\e_0,\e_0)$ with a smooth curve
$$
(-{\e_0},{\e_0}) \to X^*, \qquad s \mapsto \mu_s
$$
satisfying $\mu_0=0$.
\end{enumerate}
Here $\cO$ is defined as in (\ref{eq:def-cO}).
Since  $(\lambda_s,\varphi_s)$ is a solution to \eqref{eq:final-eq-to-solve} for every $s\in (-\e_0,\e_0)$,  the function  $u_{\phi_s} \in C^{2,\a}(\ov{\O})$  solves the  overdetermined boundary value  problem
\begin{equation*}
\left\{
  \begin{aligned}
    -\Delta_{g_{\phi_s}} u_{\phi_s} &=1 &&\qquad \textrm{ in } \quad\Omega ,\\
u_{\phi_s} &=0,\quad \partial_\nu u_{\phi_s} =-\dfrac{\lambda_s}{n}  &&\qquad \textrm{ on }\quad \partial\Omega ,
  \end{aligned}
\right.
\end{equation*}
where $\phi_s=\l_s+\vp_s$. Recalling \eqref{eq:def-diffeom-Psi_phi}, we thus find that the map $s \mapsto (\lambda_s, \vp_s)$  and the function $u:=u_{\phi_s}\circ\Psi_{\phi_s}^{-1}: \O_{\phi_s}\to \R$
have  the properties asserted in Theorem \ref{teo1}.
\QED

\section{Periodic Cheeger sets}\label{sec:Per}
In this section, we prove Corollary \ref{cor:Cheeger}. Considering the notation of Theorem~\ref{teo1}, we therefore fix $s \in (-\eps_0,\eps_0)$. Moreover, we recall that Theorem~\ref{teo1} yields a solution $u $ of the overdetermined problem
\begin{equation}\label{eq:Proe1Che}
\left\{
  \begin{aligned}
    -\Delta u&=1 &&\qquad \textrm{ in } \quad\Omega_{\phi_s},\\
u&=0,\quad \partial_\nu u=-\dfrac{\lambda_s}{n}  &&\qquad \textrm{ on }\quad \partial\Omega_{\phi_s},
  \end{aligned}
\right.
\end{equation}
where $\phi_s=\l_s+\vp_s$. In the following, we put $E_s:=\Omega_{\phi_s}$.
We need the following property which follows by a very simple application of the P-function method, see e.g. \cite{Sperb}. We include a proof for the convenience of the reader.

\begin{Lemma}
\label{eq:est-nab-us}
We have $|\nabla u|<\frac{\lambda_s}{n}\,\,$ in $\,\,E_s$.
\end{Lemma}

\proof
Consider the function
$$
P: \overline E_s \to \R,\qquad P(x):=|\nabla u(x)|^2.
$$
It is clear, by standard elliptic regularity, that $P$ is of class $C^2$. Moreover, in $E_s$ we have, since $-\Delta u  =1$,
$$
\Delta  P =2\sum^{N}_{i,j=1}\biggl(\frac{\partial^2 u}{\partial
x_i\partial
x_j}\biggl)^2  \geq {2}\sum^{N}_{i=1}\biggl(\frac{\partial^2 u}{\partial x^2_i}\biggl)^2\geq\frac{2}{N}(\Delta u)^2=\frac{2}{N}.
$$
Hence $\Delta P>0$ in $E_s$, and thus $P$ attains its maximum only on $\de E_s$ by the strong maximum principle. Since $P \equiv \frac{\lambda_s^2}{n^2}$ on $\partial E_s$ by (\ref{eq:Proe1Che}), the claim follows.
\QED

\noindent {\bf Proof of Corollary~\ref{cor:Cheeger} (completed).}
Since the domain $E_s$ is $2\pi$-periodic and symmetric in
$t_1,\dots,t_m$, the solution $u$ of  \eqref{eq:Proe1Che} is
$2\pi$-periodic and even in $t_1,\dots,t_m$. Next, let $a, b \in
\pi\Z^m$ with $a_i < b_i$ for $i=1,\dots,m$, and let $S_a^b$ be
defined as in (\ref{eq:defsalphabeta}). Then $\partial S_a^b$ can be
decomposed into a disjoint union $\partial S_\tau  = K \cup S^1 \cup
\dots \cup S^m$, where
$$
S^i:= \R^n \times \{t \in \R^m  \;:\; t_i \in \{a_ i, b_i\},\; t_j \in (a_ j , b_j) \; \text{for $j \not = i$}\}\qquad \text{for $i=1,\dots,m$,}
$$
and $K$ has zero $(N-1)$-dimensional Hausdorff measure.  By the properties of $u$ listed above, we then have
\begin{equation}
  \label{partialti-der}
\frac{\partial u}{\partial {t_i}}  \equiv 0 \qquad \text{on $S^i$ for $i=1,\dots,m$.}
\end{equation}
Next, let $A\subset E_s \cap S_a^b$ be a Lipschitz open set. Then $H^{N-1}$-almost
everywhere on $\partial A$ the outer unit normal $\nu_{\text{\tiny $A$}}$ of $A$ is well-defined, and $H^{N-1}$-almost everywhere on $\partial A \cap S^i$ it coincides with $(0,e_i)$ or $(0,-e_i)$, where $0 \in \R^n$ and $e_i$ denotes the $i$-th coordinate vector in $\R^m$. Consequently, (\ref{partialti-der}) implies that
\begin{equation}
  \label{partialti-der-1}
\partial_{\nu_{\text{\tiny $A$}}}u\equiv 0 \qquad \text{$H^{N-1}$-almost everywhere on $\partial A \cap \partial S_{a}^b$.}
\end{equation}
Since $u$ satisfies (\ref{eq:Proe1Che}), the divergence theorem and (\ref{partialti-der-1}) yield the inequality
$$
|A|=- \int_A \Delta u \,dx = - \int_{\de A}\partial_{\nu_{\text{\tiny $A$}}}u \:d\s = -
 \int_{\de A \cap S_a^b}\partial_{\nu_{\text{\tiny $A$}}} u  \:d\s  \leq
\int_{\de A \cap S_a^b} |\nabla u| d\s.
$$
Hence Lemma~\ref{eq:est-nab-us} implies that
$$
|A| \le \frac{\lambda_s}{n}H^{N-1}(\partial A\cap S_a^b)= \frac{\lambda_s}{n}P(A, S_a^b),
$$
whereas equality holds if and only if $H^{N-1}(\de A \cap S_a^b \cap E_s)=0$, i.e. if
$A= E_s \cap S_{a}^b$. This implies that $E_s$ is uniquely self-Cheeger relative to $S_a^b$ with corresponding relative Cheeger constant $h(E_s,S_{a}^b)= \frac{n}{\lambda_s}$, as claimed.  The proof of Corollary \ref{cor:Cheeger} is finished.
\QED
%
%
%
%
\vspace{0.2cm}

\textbf{Acknowledgement:}
M.M.F. is supported by the Alexander von Humboldt Foundation. Part of the paper was written while I.A.M. visited the Institute of Mathematics of the Goethe-University Frankfurt. He wishes to thank the institute for its hospitality and the German Academic Exchange Service (DAAD) for funding of the visit within the program 57060778. Also T.W. wishes to thank DAAD for funding within the program 57060778.\\

\textbf{Conflict of Interest:} The authors declare that they have no conflict of interest.

\end{document}